\documentclass[a4paper, 12pt, leqno]{article}
\usepackage{amsmath, amssymb}
\usepackage{pict2e}
\numberwithin{equation}{section}
\title{The Hele-Shaw model for the phase field method and the vortex sheet approach}
\author{Jun-ichi Koga}
\begin{document}
\maketitle
\tableofcontents
\noindent\\
\hspace{-0.2truecm}
Various thermodynamical phenomena have occurred with change of pressure and temperature, volume. We can choose these parameters but not these constraints, in order to need the thermodynamics with physical properties in the fields of various situations. The Hele-Shaw flow is one of the classical fluid dynamics problems in chemical engineering. Indeed, the Stokes flow is represented as the Stokes equations which is simplified from the famous Navier-Stokes equations. However we have often to deal with the Hele-Shaw flow with the Darcy's law. This law is derived by the Stokes equations. The basic ideas are the executive concepts. Our original work would like to solve the unknown relationship between the phase field method and the vortex sheet approach to the Hele-Shaw model.\\\\
\S 1 is devoted to Introduction.\\
In \S 2, we provide the mathematical formations for the Darcy's law, the Navier-Stokes equations, the phase field equations and the vortex sheet approach.\\
In \S 3, we calculate the two kinds of handy simulations. More precisely, the former is to discuss the Darcy's law (the Laplace equations with potential) from the phase field equations. The latter is to show the equivalence between the Hele-Shaw model and the vortex sheet approach.\\
\S 4 tells us that what models the Hele-Shaw model analogous with the Saffman-Taylor instability, explains whose equations control the Hele-Shaw model and its meaning, and we write the open questions.
\section{Introduction}
One of the most famous fellows of the Royal Society, Henry Selby Hele-Shaw proposed a certain kind of ``Stokes flow''s between two separated parallel flat plates. That gap is sufficiently small. We consider due to Caginalp [3] the Hele-Shaw model for the phase field equations. 
In physics, Saffman and Taylor\- re\-ported that fingering phenomenon for the Stokes flow which had the same gap was equivalent to the Hele-Shaw model. Mclean and Saffman [4] considered the so-called Saffman-Taylor instability with their surface tension. Pitts [5] calculated the constant $b^2/12\mu$, here $b$ is the transverse thickness of the gap and $\mu$ is the viscosity values (see Sec 2 in this paper). Tryggvason and Aref [8] numerically investigated the interfacial structure of the Hele-Shaw flow through the Stokes flow beside the two different immiscible fluids, where their boundary conditions are periodic. That is to say, their numerical experiments suggested that the Saffman-Taylor instability was analogous with Hele-Shaw cell model and also said that fingering in porous media occurs during gas injection, as used in the underground storage of gas, and in water infiltration. However their paper did not treat the dendrite growth. From the different view, they [9] discussed this instability (not the dendrite growth) with small amplitude. On the other hand, mathematically, Alikakos {\it et. al.} [1] studied the convergence of the phase field equations (Cahn-Hillard equations) to the Hele-Shaw model with the matched asymptotic expansion in 1994. Caginalp had ``had'' studied before this time that the phase field equations of ``Caginalp type'' converged to the Hele-Shaw model with the asymptotic limit (scaling limit). And one can know that the At wood ratio for viscosity tending to one the Hele-Shaw dynamics is analogous with the one-sided solidification (Stefan problems). In this paper, we reveal the unknown relationship between the vortex sheet approach and the phase field method to the Hele-Shaw model.\\\\
{\bf Theorem 1}\quad {\it Provided that the phase field method describes the Hele-Shaw model, the vortex sheet approach can also describe the Hele-Shaw model. Namely, one can derive the vortex-sheet approach from the phase field method through the Hele-Shaw model.}\\\\  
{\it Our History of Inspirations.}\quad Our origin of inspirations is at Saffman and Taylor [7]. In [7], they researched the movements of the fluids inside the porous media. Our study in this paper is based on Caginalp's work [3]. Especially, we refer to Sec. VII. in Caginalp [3].
\subsection{Discovery of Hele-Shaw Equations}
In this subsection, we introduce the Bear's way, including the Imary's average-methods.\\\\
Although Darcy originally had derived these equations, we calculate and present the Darcy's law as an empirical relationship based on experiments of steady flow in a vertical column of porous media, which is analogous with the Hele-Shaw flow. In this subsection, we have expected that we introduce the various extensions of laws ,which, for the three-dimensional Hele-Shaw flow and for the unsteady flow, etc. are heuristic generalizations. However, as one of them, Imary (1958) studied the derivations of the Darcy's law from the average of the Navier-Stokes equations.\\
Let us explain these theoretical calculations. Imary rewrote the Navier-Stokes equations as
\begin{equation}
g\frac{\partial E}{\partial x}=\frac{1}{2}\frac{\partial (V_y^2+V_z^2-V_x^2)}{\partial x}+\frac{\partial (V_xV_y)}{\partial y}+\frac{\partial (V_xV_z)}{\partial z}-\frac{\partial V_x}{\partial t}+\nu\left(\frac{\partial^2 V_x}{\partial x^2}+\frac{\partial^2V_x}{\partial y^2}+\frac{\partial^2 V_x}{\partial z^2}\right)
\end{equation}
and similar equations for $y$ and $z$ directions. The total energy of flow per unit weight of fluid is:
$$
E:=\varphi +\frac{V^2}{2g}=z+\frac{p}{\rho g}+\frac{V^2}{2g}.
$$
Imary (1958) then computed the space average for a model, made of spheres of diameters $d$. $n=L/(L+d)$. Now we also know the continuity equations, for the incompressible homogeneous fluid,
\begin{equation}
\frac{\partial V_x}{\partial x}+\frac{\partial V_y}{\partial y}+\frac{\partial V_z}{\partial z}=0
\end{equation}
where $\displaystyle\overline{V_zV_y}=0$, $\overline{V_xV_z}=0$, $\displaystyle\frac{\partial(\overline{V_xV_y})}{\partial_y}=0$ and $\displaystyle\frac{\partial(\overline{V_xV_z})}{\partial_z}=0$. The specific discharge and its time derivative are given by $\overline{V}_x=q/n$ and $\partial\overline{V}_x/\partial t=(1/n)\cdot (\partial q/\partial t)$.
Thus we get the following equations:
\begin{equation}
\frac{\partial^2\overline{V}_x}{\partial y^2}+\frac{\partial^2\overline{V}_x}{\partial z^2}=-\left[\beta \frac{(1-n)^2}{n^3d^2}\right]q
\end{equation}
where $\beta$ is a numerical shape factor that depend on the shape of the grains but not on the porosity or the diameter. Since $\mathrm{div} V=0$, we have:
\begin{equation}
\frac{1}{2}\frac{\partial (V_x^2)}{\partial y^2}=V_x\frac{\partial V_x}{\partial x}=-V_x\frac{\partial V_y}{\partial y}-V_x\frac{\partial V_z}{\partial z}
\end{equation}
His long argument yields
\begin{equation}
\overline{\jmath}_x=Wq_x
\end{equation}
which is the Darcy's law, where $\overline{\jmath}_x=-\partial\overline{E}/\partial x$.
\subsection{Interpretations}
We discuss the interpretations of the Hele-Shaw flow analogous with the Saffman-Taylor instability and the detail part of our perspective points. In this subsection, one can obtain the concept for the physical properties of a certain kind of fluid dynamics. More precisely, we analyze the moving interface between immiscible different fluids. The Darcy's laws on the Hele-Shaw problems are equivalent to the Stefan conditions on the Stefan problems. Similarly, the jump-behavior occurs on the boundary between two regions (e.g. water and ice, gricelin and air, respectively) (see e.g. Sec 2, SubSec 3.2 Eq.(3.11)).\\
Bachelor explained, without bothering the physical complicated concepts, straightly that as the field of fluid dynamics its overview shows the instabilities which fluids have in order to consult with the vortex sheet approach. Tryggvason and Aref directly and numerically calculate the phenomena by using these methods and continue the second paragraph.\\
Now the first view in their natural experiments [7] shows that fingering phenomenon occurs in the porous media. In this discussion, one can get the validity of Darcy's law. Then {\it where does the Darcy's law holds?} In physics, we interpret that this law holds when the porous media beside which the fluids move exists. We can investigate this terminology of various situations. All we have to do is to show the fundamental calculations and theoretical rationalizations (see Sec 3).\\
The second key is to determine the instability of mysterious fluids-movements. Namely, these instabilities have the wonderful interfacial structure, but almost symmetric structure (e.g. Rayleigh-Taylor instability). {\it What is an instability?} Because these phenomena depend on the linear stabilities.
\section{Mathematical Formations}
We introduce the Darcy's law (here $b$ is a gap)
\begin{equation}
{\bf u}=-\frac{b^2}{12\mu}\nabla\pi +\rho\mathbf{g},\,\,\nabla\cdot\mathbf{u}=0
\end{equation}
from the Navier-Stokes equations (See Bear [2], Saffman and Taylor [7])
\begin{equation}
\frac{\partial {\bf u}}{\partial t}+\nabla\cdot {\bf uu}=\nabla\cdot\mu(\nabla {\bf u}+\nabla^T {\bf u})+\nabla\pi +\rho\mathbf{g},\,\, \nabla\cdot\mathbf{u}=0,\,\,x\in\mathbb{R}^3.
\end{equation}
However, instead of these momentum equations, we consider the phase field equations
\begin{equation}
\left\{
\begin{array}{l}
\displaystyle u_t=\frac 12 \varphi_t+\Delta u\\\\
\displaystyle\varphi_t=\Delta\varphi +\frac 12 \left(\varphi -\varphi^3\right)+\lambda u
\end{array}
\right.
\end{equation}
to show the equivalence between the Hele-Shaw flow and, the model, which is reduced to by the phase field method, namely, without the unknown velocity functions like the Navier-Stokes equations we can deal with our investigation for the interface of two different immiscible fluids in the porous media, ``directly''.\\
As a result, we can obtain the following equations:
\begin{equation}
\begin{array}{l}
\displaystyle\Delta u=0\,\,\mathrm{in}\,\,\Omega_1,\Omega_2,\\\\
\displaystyle [\nabla u]^+_-=-\frac{\ell}{K}v\,\,\mathrm{on}\,\,\Gamma,\\\\
\displaystyle u=\left[\frac{-\alpha v-\kappa}{4}\right]\sigma\,\,\mathrm{on}\,\,\Gamma,\\\\
\end{array}
\end{equation}
subject to suitable initial and boundary conditions. Here $u$ is the pressure instead of temperature and the other physical constants have a different meaning, while $\sigma$ is the interfacial tension between the fluids. We consider due to Caginalp [3] the two cases, $\alpha >0$ and $\alpha =0$.\\
In addition, we also know that the vortex sheet approach [8] is considered as
\begin{equation}
u=\frac{1}{2\pi}\int_\Omega\frac{\widetilde{z}\times\left(x-x^\prime (s,t)\right)}{\left|x-x^\prime (s,t)\right|^2}\gamma (s,t)ds
\end{equation}
where $s$ is arclength, $\gamma (s,t)$ is the vortex-sheet strength and $\widetilde{z}$ is a unit vector perpendicular to the plane of motion.
\section{Calculations for Proof of Theorem 1}
We divide Proof of Theorem 1 into twofolds: the fundamental calculations and theoretical rationalizations. The former calculations are almost due to Caginalp's work [3] ({\bf Lemma 2}). The latter calculations are our original work almost due to Tryggvason and Aref [8,9] ({\bf Lemma 3}). More precisely, this is based on the vortex sheet approach to the phase field methods.
\subsection{Fundamental Calculations}
We discuss due to Caginalp [3] the outer- and inner- expansions.\\\\
{\bf Lemma 2.}\quad {\it The phase field equations $(2.3)$ converges to the Hele-Shaw type model $(2.4)$ as a singular limit.}\\\\
{\it Proof.}\quad Now, using the parameter $\epsilon_h=\xi^2$ we rewrite the coupled phase field equations (2.3) as (for more details, see Appendix)
\begin{equation}
\begin{array}{l}
\displaystyle\epsilon_hu_t+\frac{c^2_2}{2}\varphi_t=\Delta u\\\\
\displaystyle\alpha\epsilon_h^2\varphi_t=\epsilon_h^2\Delta\varphi +\frac{1}{2}\left(\varphi -\varphi^3\right)+2u\epsilon_h.\\\\
\end{array}
\end{equation}
we have the outer expansion: For $O(1)$,
\begin{equation}
\begin{array}{l}
\displaystyle\frac{c^2_2}{2}\varphi_t^0=\Delta u^0\\\\
\displaystyle\frac{1}{2}\left(\varphi^0 -\varphi_0^3\right)=0
\end{array}
\end{equation}
for $O(\epsilon_h^2)$,
\begin{equation}
\begin{array}{l}
\displaystyle u^0_t+\frac{c^2_2}{2}\varphi_t^1=\Delta\varphi^1\\\\
\displaystyle \left(\frac{1}{2}\left(\varphi -\varphi^3\right)\right)^\prime\varphi^1+2u^0=0.
\end{array}
\end{equation}
Notice that the following expansions hold:
\begin{equation}
\begin{array}{l}
u(x,y,t,\epsilon_1)=u^0+\epsilon_1u^1+\epsilon^2_1\cdots +u^2+\cdots\\\\
\varphi (x,y,t,\epsilon_1)=\varphi^0+\epsilon_1\varphi^1+\epsilon_1^2+\cdots\\\\
s(x,y,t,\epsilon_1)=s^0+\epsilon_1s^1+\epsilon_1^2+\cdots +\varphi^2+\cdots.
\end{array}
\end{equation}
Indeed, the $O(1)$ outer expansion has solutions $\varphi^0=\pm 1$ or $0$ and 
\begin{equation}
\Delta u^0=0 (r\not=0).
\end{equation}
Because of the assumption (or definition),
\begin{equation}
\int_{-\infty}^\infty\left\{\psi^\prime (z)\right\}^2dz\equiv\sigma
\end{equation}
where we refer to the moving coordinate system $(r,z)$ and $z:=r/\epsilon_1$. So
\begin{equation}
\displaystyle\left[u^0_r\right]^+_-=-c^2_2v^0
\end{equation}
and the first matching relations yield 
\begin{equation}
c^2_2=\frac{\ell}{K}
\end{equation}
thus we have from (3.7), (3.8)
\begin{equation}
\displaystyle\left[u^0_r\right]^+_-=-\frac{\ell}{K}v_0.
\end{equation}
As in the same manner,
\begin{equation}
\displaystyle u^0\left|_{\Gamma_\pm}\right.=\left[\frac{-\alpha v^0-\kappa^0}{4}\right]\sigma.
\end{equation}
Note that $\psi$ is determined by the physical properties. In mathematics, $\psi$ is defined by $\psi\equiv\varphi^0$. Therefore the above equations (3.5), (3.9), (3.10) rationalize the Hele-Shaw model from the phase field equations (2.3). $\Box$
\subsection{Theoretical Rationalizations}
We consider due to Tryggvason and Aref [8] the vortex sheet approach to the phase field method. Directly, we calculate the Hele-Shaw model as the vortex sheet approach. First we will derive the phase field equations from the Darcy's law and, next submit the integral formula as the vortex sheet approach.\\\\
{\bf Lemma 3.}\quad {\it Darcy's law leads us to the vortex sheet approach and this approach to the the Hele-Shaw model. Inversely, we can obtain the Darcy's law from the vortex-sheet approach.}\\\\
{\it Proof.}\quad Then due to Tryggvason and Aref [8], The difference of gradient of pressure which rewrites the Darcy's law (2.1),
\begin{equation}
\left(\nabla p_2-\nabla p_1\right)\cdot\hat{s}=-\left(\frac{12\mu_2}{b^2}\mathbf{u}_2-\frac{12\mu_1}{b^2}\mathbf{u}_1\right)\cdot\hat{s}-\mathbf{g}\left(\rho_2-\rho_1\right)\cdot\hat{\jmath}\cdot\hat{s}
\end{equation}
which correspond to the Stefan conditions, on the Hele-Shaw model, equivalent to the phase field equations, instead of temperature where
\begin{equation}
\gamma =\frac{\Delta\mu}{\mu}\mathbf{U}\cdot\hat{s}+\frac{b^2g}{12\overline{\mu}}\Delta\rho\mathbf{g}\cdot\hat{s}+\frac{b^2}{12\overline{\mu}}\nabla (\Delta\rho)\cdot\hat{s}.
\end{equation}
for the viscosity of immiscible different fluids, by the way, note that the stream functions $\widetilde{\psi}$ corresponding to the velocity field equation (2.5) satisfies a Poisson equation
\begin{equation}
\Delta\widetilde{\psi}=-\omega.
\end{equation}
For small spacing difference between vortices $\Delta s$, we define $\Gamma_i$
\begin{equation}
\Gamma_i=\gamma (s_i)\Delta s.
\end{equation}
Then the $z$-direction of vorticity is
\begin{equation}
\omega (x,t)=\sum_i\Gamma_i\delta(x-x_i(t))
\end{equation}
where $x_i$ is the position vector of elemental vortex $i$. In order to move vortex $i$ according to
\begin{equation}
\frac{dx_i}{dt}=\mathbf{U}(s_i).
\end{equation}
Then
\begin{equation}
\gamma=\frac{\Delta\mu}{\overline{\mu}}\mathbf{U}\cdot\hat{s}+\left(\frac{\Delta\mu}{\overline{\mu}}\mathbf{U}_\infty +\frac{\Delta\rho g b^2}{12\overline{\mu}}\right)\hat{\jmath}\cdot\hat{s}+\frac{ab^2}{12\overline{\mu}}\frac{\partial}{\partial s}\left(\frac{1}{R_\parallel}\right).
\end{equation}
By neglecting the viscous stresses at moving interface in a real Hele-Shaw cell, we can rewrite
\begin{equation}
p_2-p_1=\alpha\left(\frac{1}{R_\parallel}+\frac{1}{R_\perp}\right),
\end{equation}
here $\alpha$ is the surface tension coefficient, $R_\parallel$ is the radius of curvature of the interface in the plane of motion and $R_\perp$ is the radius of curvature in the direction perpendicular to the parallel plates. Mclean and Saffman [4] proposed that under the assumption $R_\parallel$ is constant, we simply obtain
\begin{equation}
\nabla (\Delta p)\cdot\hat{s}=\alpha\frac{\partial}{\partial s}\left(\frac{1}{R_\parallel}\right).
\end{equation}
Then the speed of fluids is defined by
$$
U_*:=\left|\frac{(\mu_2-\mu_1)U_\infty+\frac{1}{12}(\rho_2-\rho_1)gb^2}{\mu_2+\mu_1}\right|.
$$
here $W$ is the width of the Hele-Shaw cell.
Therefore, we can get
\begin{equation}
\widetilde{\gamma}=2A(\widetilde{U}\cdot\hat{s})\pm2\hat{\jmath}\cdot\hat{s}+2B\frac{\partial}{\partial s}\left(\frac{1}{R_\parallel}\right)
\end{equation}
as the non-dimensional form of the vortex-sheet strength $\widetilde{\gamma}$, where $A$ is $\frac{\mu_2-\mu_1}{\mu_2+\mu_1}$, $B$ is $\frac{ab^2}{12U_*W^2\overline{\mu}}$.\\
Now notice that one can consider the incompressible fluid flow with $u(x,t)$ and the pressure $\pi(x,t)$, then the incompressible Navier-Stokes equations are written as
\begin{equation}
u_t+(u\cdot\nabla)u=-\nabla \pi+\nu\Delta u,\,\,\nabla\cdot u=0,
\end{equation}
where $\nu$ is the fluid viscosity. As $\nu\to 0$, and taking the curl in this limit yield the vorticity equations
\begin{equation}
\omega_t+(u\cdot\nabla)\omega=(\omega\cdot\nabla)u
\end{equation}
which and (3.20) arrive at the equation (2.5).\\\\
Inversely, we will derive the Darcy's law from the integral formula (2.5) as the vortex sheet approach. At first, we calculate the vortex sheet length $\gamma$,
\begin{equation}
\gamma=(u_1-u_2)\cdot\hat{s}.
\end{equation}
Note that Birkhoff's integral formula
\begin{equation}
\mathbf{U}(s,t)=\frac{1}{2\pi}\mathrm{P.V.}\int\frac{\widetilde{z}\times (x(s,t)-x(s^\prime,t))}{|x(s,t)-x(s^\prime,t)|^2}\gamma(s^\prime,t)ds^\prime,
\end{equation}
which geometrically allows us to estimate and calculate (3.17), (3.20), (3.23) with this integral formula (see figure 4.1 Model), more precisely, the equation (3.13) makes the vorticity equation, the equation (3.15) creates the (3.24) and these processes provide us with the equation (3.19), that, in addition, makes the equation (3.18). Thus this equality guarantees the most left hand-side of the equation (3.11). Naturally, we can obtain the generalized Darcy's law (2.1).$\Box$\\\\
{\it General Remark.}\quad We recall that the two following relations are formally equivalent: The difference of gradient of the Darcy's laws,
\begin{equation}
\left(\nabla p_2-\nabla p_1\right)\cdot\hat{s}=-\left(\frac{12\mu_2}{b^2}\mathbf{u}_2-\frac{12\mu_1}{b^2}\mathbf{u}_1\right)\cdot\hat{s}-\mathbf{g}\left(\rho_2-\rho_1\right)\cdot\hat{\jmath}\cdot\hat{s},
\end{equation}
and the Stefan conditions,
\begin{equation}
\ell v=K(\nabla u_S-\nabla u_L)\cdot\widehat{\mathbf{n}}, x\in\Gamma(t).
\end{equation}
\section{Main Result}
\subsection{Comments and Models}
{\it Comments.} It is well-known that the Hele-Shaw model is analogous with the Saffman-Taylor instability. Their analogy suggests that the governing equations are the same as the other at the same time (see SubSec 4.2). As in the model, we consider water or oil inside the porous rock or sand, respectively. In physics, we can observe this phenomenon where filtrations of oil is removed. On the other hand, the vortex sheet approach is suitable for the investigation of these instabilities due to Tryggvason and Aref [8]. However the phase field equations are developed as the novel calculation method to numerically and theoretically simulate the interfacial structure of two immiscible fluids. The following figure represents the vortex sheet approach. For more details, there are two separated parallel flat glass-like plates and through so small gap that the Stokes flow occurs with the fluids.\\\\
{\it Models.}\quad We study the Hele-Shaw flow between so narrow parallel flat plates that the fluid is governed as the Darcy equations in $2+1$-dimension. Although these approximations do not seem correct, it becomes realistic in this field of fluid dynamics (see Pitts [4]). Because the gap is sufficiently small. Our treatment is perspective, and the vortex sheet approach is similar to the phase field method to the Hele-Shaw model. In this discussion, we solve this constraint to prove this terminology. Indeed, the vortex-sheet approach comes from the jump-behavior of velocities.
$$
\Gamma_i=\{v\}\cdot ds
$$
where $\Gamma_i$ is a sum of vortex filaments on the area surrounded by $C_i$, $\{v\}$ is the jump of the velocities. For this equation, we apply the Darcy's law,
$$
\Gamma_i^\prime=\left(\nabla \pi_2-\nabla \pi_1\right)\cdot\hat{s}
$$
which is equivalent to the Stefan conditions of the phase field method as a singular limit, where the viscosities are averaged. The derivation of the degenerate heat equations is not the most important problem. The ``Stefan conditions'' are the most important calculations in the realistic fields. Proof of Lemma 1 and proof of Lemma 2 yield proof of Theorem 1 from these perspective views. 
\\\\\\\\\\
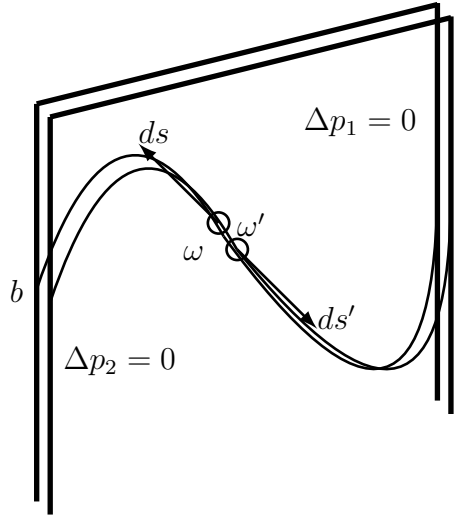
\begin{center}
\begin{picture}(0,-500)
\linethickness{2pt}
\put(-60,0){\line(40,10){150}}
\put(-55,-5){\line(40,10){150}}
\put(-60,0){\line(0,-1){150}}
\put(-55,-5){\line(0,-1){150}}
\put(90,38){\line(0,-1){150}}
\put(95,33){\line(0,-1){150}}
\linethickness{1pt}
\put(-70,-75){$b$}
\qbezier(-60,-70)(-30,20)(10, -50)
\qbezier(10,-50)(90,-160)(90, -28)
\qbezier(-55,-75)(-25,15)(15, -55)
\qbezier(15,-55)(95,-155)(95, -33)
\put(15,-55){\circle{8}}
\put(8,-45){\circle{8}}
\put(15,-55){\vector(2,-2){30}}
\put(8,-45){\vector(-1,1){30}}
\put(45,-85){$ds^\prime$}
\put(-22,-15){$ds$}
\put(15,-50){$\omega^\prime$}
\put(-5,-59){$\omega$}
\put(-50,-100){$\Delta p_2=0$}
\put(40,-10){$\Delta p_1=0$}
\end{picture} 
\end{center}
\vspace{5truecm}
\subsection{Governing Equations and Its Meaning}
The governing equations are three patterns: The Darcy's law
\begin{equation}
\mathbf{u}=-\frac{b^2}{12\mu}\nabla\pi,\,\,\nabla\cdot\mathbf{u}=0
\end{equation}
and its meaning tells us that the Stokes equations approximate the fluid flow as linearly.\\\\
{\it Remark.}\quad The same law holds for the flow beside the porous media and the Hele-Shaw flow at the same time.\\\\
: The phase field equations
\begin{equation}
\left\{
\begin{array}{l}
\displaystyle u_t=\frac 12 \varphi_t+\Delta u\\\\
\displaystyle\varphi_t=\Delta\varphi +\frac 12 \left(\varphi -\varphi^3\right)+\lambda u
\end{array}
\right.
\end{equation}
where $u$ is the pressure instead of temperature, $\varphi$ is the order parameters which control the interface between two regions.\\
:The vortex sheet approach
\begin{equation}
u=\frac{1}{2\pi}\int_\Omega\frac{\widetilde{z}\times\left(x-x^\prime (s,t)\right)}{\left|x-x^\prime (s,t)\right|^2}\gamma (s,t)ds
\end{equation}
where $s$ is arclength, $\gamma$ is the vortex-sheet length, and $\widetilde{z}$ is a unit vector perpendicular to the plane of motion. 
\subsection{Open Questions}
We have the following questions:
\begin{itemize}
\item What is the equivalence between the vortex sheet approach and the phase field method to the Hele-Shaw model?
\item What is the difference between the Hele-Shaw flow and the Saffman-Taylor instability?
\end{itemize}
\appendix
\section{Appendix}
In the equation (2.3), we make the following parameters
\begin{equation}
\epsilon^2\equiv\xi^2 a,\,\,\alpha=\mathrm{fixed},\,\,\xi, a\to 0,\,\, \rho\equiv r/\epsilon
\end{equation}
and define the following expansions using $\displaystyle\frac{1}{2}\left(\phi -\phi^3\right)$.
\subsection{Definition-Inner expansions}
\begin{equation}
-\alpha v\epsilon\phi_\rho\cong\phi_{\rho\rho}+\epsilon\kappa\phi_\rho +\cdots (O(\epsilon^2))+\frac{1}{2}\left(\phi -\phi^3\right)+2au.
\end{equation}
to rewrite the coupled phase field equations.\\
Then we write the inner expansions:
\begin{equation}
\begin{array}{ll}
u(x,y,t,\epsilon_1) & =U(z,s,t,\epsilon_1)\\
                      &=U^0(z,s,t)+\epsilon_1U^1(z,s,t)+\cdots\\\\
\varphi (x,y,t,\epsilon_1)&=\phi(z,s,t,\epsilon_1)\\
                                &=\phi^0(z,s,t)+\epsilon_1U^1(z,s,t)+\cdots
\end{array}
\end{equation}
\subsection{Definition-Outer expansion}
Before defining the outer expansions, we appeal the following relations:
\begin{equation}
\begin{array}{l}
u(x,y,t,\epsilon_1)=u^0+\epsilon_1u^1+\epsilon^2_1\cdots +u^2+\cdots\\\\
\varphi (x,y,t,\epsilon_1)=\varphi^0+\epsilon_1\varphi^1+\epsilon_1^2+\cdots\\\\
s(x,y,t,\epsilon_1)=s^0+\epsilon_1s^1+\epsilon_1^2+\cdots +\varphi^2+\cdots,
\end{array}
\end{equation}
which is determined by the mathematical calculations.
\subsection{Solutions of ODE}
Now let us calculate the following ordinary differential equations
\begin{equation}
\phi^0_{\rho\rho}+\frac{1}{2}[\phi^0-(\phi^0)^3]=0
\end{equation}
satisfying
\begin{equation}
\phi_0(\rho)=\tanh(\rho/2). 
\end{equation}
Multiplying both sides (A.5) by $\phi_\rho$ and integrating in $\rho$,
\begin{equation}
\frac{1}{2}\left(\frac{d\phi}{d\rho}\right)^2=\frac{1}{8}\left(1-\phi^2\right)^2.
\end{equation}
So
\begin{equation}
\frac{d\phi}{d\rho}=\frac{1}{2}\left(1-\phi^2\right),
\end{equation}
then we have
\begin{equation}
-\ln(1-\phi)+\ln(1+\phi)=\rho +C_1
\end{equation}
\begin{equation}
\ln\frac{1+\phi}{1-\phi}=\rho\Leftrightarrow \frac{1+\phi}{1-\phi}=e^\rho.
\end{equation}
Thus we get the desired result (A.6).

\end{document}